\documentclass[11pt]{article}      % onecolumn (second 
\usepackage{graphicx}
\usepackage{latexsym}
\usepackage{amssymb}
\usepackage[mathscr]{eucal}
\usepackage{color}
\newcommand{\Ep}{\hfill\hfill{\qed}\end{proof}}
\newcommand{\Bp}{\begin{proof}}

\newcommand{\plq}{\rule{0.in}{.25in}\left[}
\newcommand{\prq}{\rule{0.in}{.25in}\right]}

\newcommand{\Div}{\mbox{\rm div}\,}

\newcommand{\supp}{\mbox{\rm supp}\,}

\newcommand{\curl}{\mbox{\rm curl}\,}

\newcommand{\Int}[2]{{\displaystyle \int_{ #1}^{ #2}}}

\newcommand{\Sup}[1]{{\displaystyle \sup_{#1}}}

\newcommand{\Sum}[2]{{\displaystyle \sum_{#1}^{#2}}}

\newcommand{\beea}{\begin{eqnarray}}

\newcommand{\eeea}{\end{eqnarray}}

\newcommand{\ms}{\medskip\smallskip}

\newcommand{\ds}{\displaystyle}

\newcommand{\bfe}{{\mbox{\boldmath $e$}} }

\newcommand{\bfz}{{\mbox{\boldmath $z$}} }

\newcommand{\0}{{\mbox{\boldmath $0$}} }

\newcommand{\BF}{\begin{footnotesize}}

\newcommand{\EF}{\end{footnotesize}}

\newcommand{\pde}[2]{{\displaystyle \frac{\mbox{$\partial #1$}}{\mbox{$\partial #2$}}}}

\newcommand{\ode}[2]{{\displaystyle \frac{\mbox{$d #1$}}{\mbox{$d #2$}}}}

\newcommand{\bi}{\begin{itemize}}

\newcommand{\ei}{\end{itemize}}

\newcommand{\ed}{\end{document}}

\newcommand{\be}{\begin{equation}}

\newcommand{\ba}{\begin{array}}

\newcommand{\ea}{\end{array}}

\newcommand{\ee}{\end{equation}}

\newcommand{\eeq}[1]{\label{eq:#1}\end{equation}}

\newcommand{\real}{{\mathbb R}}%\newcommand{\F}{{\rm I\!\!\,F}}

\newcommand{\nat}{{\mathbb N}}

\newcommand{\bfpsi}{\mbox{\boldmath $\psi$}}

\newcommand{\bfomega}{\mbox{\boldmath $\omega$}}

\newcommand{\bfxi}{\mbox{\boldmath $\xi$}}

\newcommand{\bfx}{\mbox{\boldmath $x$}}

\newcommand{\bfy}{\mbox{\boldmath $y$}}

\newcommand{\bfv}{{\mbox{\boldmath $v$}} }

\newcommand{\bfu}{{\mbox{\boldmath $u$}} }

\newcommand{\bfw}{{\mbox{\boldmath $w$}} }

\newcommand{\bff}{{\mbox{\boldmath $f$}} }

\newcommand{\bfc}{{\mbox{\boldmath $c$}} }

\newcommand{\bfQ}{{\mbox{\boldmath $Q$}} }

\newcommand{\bfH}{{\mbox{\boldmath $H$}} }

\newcommand{\bfh}{{\mbox{\boldmath $h$}} }

\newcommand{\cala}{{\cal A}}

\newcommand{\cals}{{\cal S}}

\newcommand{\bfV}{{\mbox{\boldmath $V$}} }

\newcommand{\bfF}{{\mbox{\boldmath $F$}} }

\newcommand{\bfb}{{\mbox{\boldmath $b$}} }

\newcommand{\bfg}{{\mbox{\boldmath $g$}} }

\newcommand{\bfn}{{\mbox{\boldmath $n$}} }

\newcommand{\half}{\mbox{$\frac{1}{2}$}}

\def\Bbb R{\real}

\def\tilde{\widetilde}

\def\bar{\overline}

\newcommand{\bfcalb}{\mbox{\boldmath ${\cal B}$}}

\newcommand{\bfcalg}{\mbox{\boldmath ${\cal G}$}}

\newcommand{\bfcalh}{\mbox{\boldmath ${\cal H}$}}

\newcommand{\bfcalf}{\mbox{\boldmath ${\cal F}$}}

\newcommand{\ED}{\end{description}}

\newcommand{\Footnote}{~\footnote}

\newcommand{\Br}{\begin{rem}\begin{rm}}

\newcommand{\Er}{\end{rm}\end{remark}}

\newtheorem{lemm}{Lemma}[section]

\newtheorem{theo}{Theorem}[section]

\newtheorem{rem}{Remark}[section]

\newtheorem{coro}{Corollary}[section]

\newtheorem{exe}{\footnotesize{Exercise}}[section]

\newcommand{\Be}{\begin{exe}\begin{footnotesize}\begin{rm}}

\newcommand{\EE}[1]{\end{rm}\end{footnotesize}\label{exe:#1}\end{exe}}

\newcommand{\Bt}{\begin{theo}\begin{sl}}

\newcommand{\Et}{\end{sl}\end{theorem}}

\newcommand{\Bl}{\begin{lemm}\begin{sl}}

\newcommand{\El}{\end{sl}\end{lemma}}

\newcommand{\eqref}[1]{{\rm (\ref{eq:#1})}}

\newcommand{\Bc}{\begin{coro}\begin{sl}}

\newcommand{\Ec}{\end{sl}\end{coro}}

\newcommand{\ET}[1]{\end{sl}\label{theo:#1}\end{theo}}

\newcommand{\EL}[1]{\end{sl}\label{lemm:#1}\end{lemm}}

\newcommand{\theoref}[1]{{\rm Theorem \ref{theo:#1}}}

\newcommand{\ER}[1]{\end{rm}\label{rem:#1}\end{rem}}

\newcommand{\EC}[1]{\end{sl}\label{coro:#1}\end{coro}}

\newcommand{\lemmref}[1]{{\rm Lemma \ref{lemm:#1}}}

\newcommand{\essup}[1]{{\rm ess}\,{{\displaystyle \sup_{\hspace*{-6mm}{#1}}}}\,}

\begin{document}
\title{Existence, Uniqueness and Asymptotic Behavior\\ of Regular Time-Periodic Viscous Flow\\ around a Moving Body: Rotational Case}
\author{Giovanni P Galdi \smallskip \\ { \small Department of Mechanical Engineering and Materials Science}\\ { \small University of Pittsburgh, USA}}
\date{}
\maketitle
\begin{abstract}
We show existence and uniqueness for small data of regular time-periodic solutions to the Navier-Stokes problem in the exterior of a rigid body, $\mathscr B$, that moves by time-periodic translational motion of the same period along a constant direction, $\bfe_1$, and spins with constant angular velocity $\bfomega$  parallel to $\bfe_1$. We also study the spatial asymptotic behavior of such solutions and show, in particular, that if $\mathscr B$ has a  net motion characterized by a non-zero average translational velocity  $\bar{\bfxi}$, then the solution exhibit a wake-like behavior in the direction $-\bar{\bfxi}$ entirely analogous to that of a  steady-state flow around a body that moves with  velocity $\bar{\bfxi}$ and angular velocity $\bfomega$.
\end{abstract}
\renewcommand{\theequation}{{0}.\arabic{equation}}
\section*{Introduction}
A bounded rigid  body, $\mathscr B$, moves in a quiescent viscous liquid, $\mathscr L$,  that fills the entire space outside $\mathscr B$. Assume that, when referred to a body-fixed frame, the motion of $\mathscr B$ is time-periodic, namely, the characteristic vectors $\bfxi$ and $\bfomega$ representing the velocity of the center of mass and  angular velocity of $\mathscr B$, respectively, are prescribed,   sufficiently smooth time-periodic functions of period $T$ ({\em $T$-periodic}). The natural question that arises is whether the liquid will also execute a $T$-periodic motion that is regular and unique, at least when the magnitude of both $\bfxi$ and $\bfomega$ is  ``small" enough.  In mathematical terms, this question amounts to find $T$-periodic (regular, unique) solutions $(\bfu,p)$ to the following system of equations
\be\ba{cc}\smallskip\left.\ba{ll}\medskip
\bfu_t-(\bfxi+\bfomega\times\bfx)\cdot\nabla\bfu+\bfomega\times\bfu+\bfu\cdot\nabla\bfu=\Delta\bfu-\nabla {p}+\bfb\\
\Div\bfu=0\ea\right\}\ \ \mbox{in $\Omega\times (-\infty,\infty)$}\\
\bfu=\bfxi+\bfomega\times\bfx\,,\ \ \mbox{at\ $\partial\Omega\times (-\infty,\infty)$}\,,
\ea
\eeq{0.1}
where $\bfu$ and $p$ are velocity and pressure fields of $\mathscr L$,   and $\Omega$ is the complement of a connected compact set  of $\real^3$ (the body $\mathscr B$). Moreover, for the sake of generality, we include also a (prescribed) $T$-periodic body force $\bfb=\bfb(x,t)$ acting on $\mathscr L$.
\par
The investigation of this problem in the ``non-trivial" case, namely, when at least one of the two characteristic vectors is not identically zero, has begun   only in the recent past with the paper \cite{GS1}.\footnote{In the case when $\mathscr B$ is motionless $(\bfxi\equiv\bfomega\equiv\0$), there is a vast literature. In particular, we refer the reader to \cite{Ma,MaPa,KoNa,Y,GaSo,KMT} and the references therein.} There, it is shown that, provided only  $\bfxi$,  $\bfomega$, and $\bfb$ have a mild degree of regularity, \eqref{0.1} has at least one $T$-periodic  weak solution a la Leray-Hopf that, in fact, is also strong in the sense of Ladyzhenskaya whenever the data possess more regularity, and their magnitude is appropriately restricted. 
\par
Nonetheless, there were two fundamental issues left open in \cite{GS1}:   {\em uniqueness} and, the somehow related problem,  {\em asymptotic spatial behavior} of solutions.  
These and related questions were successively analyzed and successfully addressed by several authors by  entirely different methods than those employed in \cite{GS1}, but, however, under more restrictive  assumptions on $\bfxi$ and $\bfomega$,  as we are about to detail. Specifically, in \cite{Ga,K1,GaKy,GaKy1,EitKye,EitKye1}, one supposes that $\bfomega\equiv\0$,  $\bar{\bfxi}\neq\0$ (the bar denoting  average over a period), and $\sup_t|\bfxi(t)|$, $\sup_t|\bfxi(t)-\bar{\bfxi}|$ ``small" enough. The restriction $\bfomega\equiv\0$ is removed in \cite{Eiter}, and replaced by $\bfomega={\bf const.}$, but the period of the  driving mechanism ($\bfxi$ and/or $\bfb$) should be   equal to $2\pi\kappa/|\bfomega|$, for some non-zero $\kappa\in \mathbb Q$. Moreover, $\bfxi(t)$ and $\bfomega$ must be parallel at all times and each one ``small" enough. These results are obtained by means of a maximal regularity approach to time-periodic problems initiated in \cite{Ga} and developed and generalized to an elegant abstract form in \cite{K2,KS}. Another different approach, based on sharp $L^p-L^q$ estimates,  was taken in \cite{NTH,HG,MH,GH}. The success of this approach, however, requires   $\bfxi$ and $\bfomega$  ``small", parallel and both constant in time but. However, unlike \cite{Eiter},  no restriction is imposed on the period $T$. Actually, this method applies also to  the more general almost-periodic case \cite{HG}.  
\par
It should be emphasized that the above assumption of $\bfxi$ being constant in time or, more generally, having a non-zero average is very restrictive from the physical viewpoint, in that it excludes the simplest time-periodic motion of $\mathscr B$, namely, an oscillation between two fixed configurations (like in a pendulum, for instance). Motivated by this observation in \cite{GARMA} and, successively, in \cite{GaLN} we have used a yet different approach to prove existence and uniqueness of regular solutions to \eqref{0.1} without any restriction on the function $\bfxi=\bfxi(t)$, other than its smoothness and ``smallness", while keeping $\bfomega\equiv\0$. Moreover, we have shown  that at large distance from $\mathscr B$  the flow field  has a distinctive steady-state profile, which, in particular, has  allowed  us in \cite{GARMA} to give a rigorous interpretation of the steady streaming phenomenon \cite{Ri}. Our method relies upon the well-posedness of the corresponding linear problem combined with a contraction-mapping argument and consists of two steps. At first, one uses the classical Galerkin method endowed with a number of ``high order energy" estimates to show the existence of a regular $T$-periodic solution to the relevant linear problem in correspondence of sufficiently smooth data. Successively, if the data decay  at large distances in a well-specified fashion, then one shows that the solutions must decay at a distinctive rate as well.        
\par
The main objective of the present  paper is to suitably extend the method introduced in \cite{GARMA, GaLN}  to cover also the case $\bfomega\not\equiv\0$. Such an extension is by no means incremental since, as is well known,  spinning of $\mathscr B$ introduces a substantial difficulty to the problem, represented by an extra term in the linear momentum equation whose coefficient may grow unbounded at large spatial distances. As a matter of fact, we are only able to treat the case when $\bfomega$ is constant and, as in all papers mentioned above, constantly parallel to $\bfxi(t)$. In such a case, under suitable assumptions of smoothness and ``smallness",  we  prove existence and uniqueness of regular solutions to \eqref{0.1} and, in addition, provide a sharp asymptotic spatial behavior of the flow field; see \theoref{3.1}.  Concerning the latter, we prove the following. Take $\bar{\bfxi}=\lambda\,\bfe_1$, with $\bfe_1$ unit vector in the direction $x_1$ and $\lambda\ge 0$. Then $\bfu(x,t)$ decays like $|x|^{-1}[1+\lambda\,(|x|+x_1)]^{-1}$, uniformly in $t\in[0,T]$. Thus, if $\lambda=0$ we find a  result similar to that shown in \cite{GARMA}, but with $\bfomega\neq\0$. However, if $\lambda>0$, the velocity field presents a ``wake" behavior in the direction opposite to $\bar{\bfxi}$,  entirely analogous to that of a  steady-state flow around a body that moves with constant velocity $\bar{\bfxi}$ parallel to its angular velocity $\bfomega$ \cite[Section XI.6]{Gab}. 
\par
We end this introductory section by pointing out two  avenues of further research. The first concerns  the weakening of the above assumptions on $\bfxi$ and $\bfomega$. Perhaps, the most intriguing question is whether the hypothesis of $\bfomega$ and $\bfxi$ being parallel at all times  is indeed necessary to obtain the well-posedness of the problem, for small data of course. Notice that, as shown in \cite{GS1}, such a request is indeed not needed for existence,   even in the case  of strong solutions. The other research address  that appears to be of some appeal especially for its  geophysical applications, is the study of analogous questions in the presence of thermal effects. In this situation, it is relevant to take $\bfxi\equiv\0$, $\bfomega=\textbf{const.}$ and a $T$-periodic distribution of temperature on the surface of $\mathscr B$. Similar   studies on simplified models  and in the Boussinesq approximation, have been performed by several authors; see \cite{His} and the references therein. It would be interesting to investigate these problems for the full system \eqref{0.1} in combination with heat-conduction in the Boussinesq approximation, or the even   more realistic pressure-dependent model recently  proposed in \cite{DePa}, \cite{PR}.   
\par
The outline of the paper is as follows. After recalling in Section 1 some known preparatory results, in the following Section 2 we prove existence, uniqueness and continuous data dependence, in a suitable function class, of solutions to the linear problem obtained by neglecting the nonlinear terms in \eqref{0.1}; see \theoref{1.1}. The challenging part to get such a result consists in proving local estimates for the $L^2$-norm of $\bfu_t$, uniformly in time. While the proof of this property is relatively straightforward when $\bfomega\equiv\0$ \cite{GaLN}, it becomes more involved if $\bfomega\not\equiv\0$; see \lemmref{Ar1} and \lemmref{1.3}. In fact, this is the only point were we need $\bfomega$  constant in time. If a similar estimate could be shown also when $\bfomega$ is time-dependent, the proof of \theoref{1.1} could  be carried out in exactly the same way, providing the same results under this more general assumption. In the final Section 3,  \theoref{3.1}, we combine \theoref{1.1} with a classical contraction mapping argument and extend the linear findings to the full nonlinear case under suitable restriction on the magnitude of the data.

\renewcommand{\theequation}{\arabic{section}.\arabic{equation}}\setcounter{equation}{0}
\section{Preliminaries}
Throughout, $\Omega$ denotes  the complement of the closure of a bounded  domain $\Omega_0\subset\mathbb R^3$,  of class $C^2$. We  take the origin of the coordinate system in the interior of $\Omega_0$.   For $R\ge R_*:=2{\rm diam}\,(\Omega_0)$, we set 
$
\Omega_R=\Omega\cap \{|x|<R\}\,,\ \ \Omega^R=\Omega\cap\bar{\{|x|>R\}}$.
If $A\subseteq \real^3$ is a domain, by $L^q (A)$,  $1\leq q \leq \infty,$  
$W^{m,q}({A}),$ $W_0^{m,q}(A)$, $m \geq 0,$  $(W^{0,q}\equiv W^{0,q}_0\equiv L^q$), we indicate usual Lebesgue and Sobolev spaces, with norms $\|.\|_{q,A}$ and $\|.\|_{m,q,A}$, respectively.\Footnote{We shall use the same font style to denote scalar, vector and tensor
function spaces.} By $P$ we denote  the (Helmholtz) projector from $L^2(A)$ onto its subspace of solenoidal (vector) function with vanishing normal component, in distributional sense, at $\partial A$.   
 We also set $\int_{A}u\cdot v=( u,v)_{A}$.  
$D^{m,2}(A)$ is the space of (equivalence classes of) functions $u$ with seminorm 
$ 
\sum_{|k|=m}\|D^k u\|_{2,A}<\infty\,.
$  By $D_0^{1,2}(A)$ we denote the completion of $C_0^\infty(A)$ in the norm $\|\nabla(\cdot)\|_2$. In the above notation,  the subscript ``$A$" will be omitted, unless confusion arises. A function $u:A\times \real\mapsto \real^3$ is 
{\em $T$-periodic}, $T>0$, if $u(\cdot,t+T)=u(\cdot\,t)$, for a.a. $t\in \real$,
 and we set
$
{\bar u}:=\frac{1}{T}\int_{0}^{T}u(t){\rm d}t\,.
$
Let $B$ be a function space endowed with seminorm $\|\cdot\|_B$, $r=[1,\infty]$, and $T>0$. $L^r(0,T;B)$ is the class of functions
$u:(0,T)\rightarrow B$ such that 
$$
\|u\|_{L^r(B)}\equiv\left\{\ba{ll}\smallskip\big( \Int{0}{T}\|u(t)\|_B^r \big)^{\frac 1r}<\infty, \ \ \mbox{if 
$r\in [1,\infty)\,;$}\\   
\essup{t\in[0,T]}\,\|u(t)\|_B <\infty, \ \ \mbox{if $r=\infty.$}
\ea\right.
$$
We also define
$$%\ba{ll}\medskip
W^{m,r}(0,T;B)=\Big\{u\in L^{r}(0,T;B): \textcolor{black}{\partial_t^ku\in L^{r}(0,T;B), \, k=1,\ldots,m}\Big\}\,.
%\ea
$$
We shall simply write $L^r(B)$ for $L^r(0,T;B)$, etc. unless otherwise stated 
Finally, for $A:=\Omega,\real^3$,  $m\ge 1$, and $\lambda\ge 0$ we set
$$\ba{ll}\medskip
[\!]f[\!]_{m,\lambda,A}:=\Sup{x\in A}\,|(1+|x|)^{m}(1+2\lambda \,s(x))^m f(x)|\,,\\ 

[\!] f [\!]_{\infty,m,\lambda,A}:=\Sup{(x,t)\in A\times (0,\infty)}\,|(1+|x|)^{m}(1+2\lambda \,s(x))^m f(x,t)|\,.\ea
$$
where $s(x)=|x|+x_1$,  $x\in\real^3$, and  the subscript $A$ will be omitted, unless necessary.
\par
Finally, we set
$$
(\cdot)_t:=\pde{(\cdot)}t\,,\ \ \partial_1(\cdot):=\pde{(\cdot)}{x_1}\,,
$$
and denote by $c$ a generic positive constant whose specific value is irrelevant and may change  even in the same line.\smallskip\par
We next collect some preliminary results whose proof can be found in the literature. We begin with the following one,  a special case of \cite[Lemma II.6.4]{Gab}
\Bl
There exists a  function $\psi_{R} \in C_0^{\infty}(\real^n)$ defined for all $R>0$ such that $0\le \psi_{R}(x)\le 1$, $x\in \real^n$, and satisfying the following properties 
$$
\displaystyle \lim_{R\to\infty}\psi_{R}(x)=1\,,\  \mbox{uniformly pointwise\,;}\ \ \
\left\|\partial_1{\psi_{R}}\right\|_{\frac32}\le C_1\,,
$$
where $C_1$  is independent of  $R.$ Moreover, the support of $\partial\psi_{R}/\partial x_j$, $j=1,\ldots,n$, is contained in $\Omega^{\frac{R}{\sqrt2}}$ and
$$
\|u\,|\nabla\psi_{R}|\,\|_{2}\le C_2\,\|\nabla u\|_{2,
\Omega^{\frac{R}{\sqrt{2}}}}\,,\ \ \mbox{for all $u\in D_0^{1,2}(\Omega)$}\,.
$$
where  $C_2$ is independent of $R$. 
\EL{1.1}
The following result can be found  in \cite[Theorem III.3.1 and Exercise III.3.7]{Gab}.
\Bl Let $\cala$ be a bounded Lipschitz domain in $\real^3$, and let $f\in L^2(\cala)$ with $\int_{\cala}f=0$. Then the problem
\be
\Div\bfz=f\,\ \mbox{in $\cala$}\,, \ \bfz\in W_0^{1,2}(\cala)\,,\ \ \|\bfz\|_{1,2}\le C_0\,\|f\|_{2}\,,
\eeq{Bog}
for some $C_0=C_0(\cala)>0$ has at least one solution. Moreover, if also $f\in W_0^{1,2}(\cala)$, then $\bfz\in W_0^{2,2}(\cala)$ and $\|\bfz\|_{2,2}\le C_0\,\|f\|_{1,2}$.
Finally, if $f=f(t)$ with $ f_t\in L^\infty(L^2(\cala))$, then we have in addition $\bfz_t\in L^\infty(W_0^{1,2}(\cala))$ and
$$
\|\bfz_t\|_{1,2}\le C_0\,\|f_t\|_2\,.
$$ 
\EL{1.1_1}
The next result is proved in  \cite[Lemma 2.2]{GS1}. 
\Bl
Let $\bfxi \in W^{\textcolor{black}{2},2}(0,T)$ be $T$-periodic, and let $\bfomega\in\real^3$. Given $\varepsilon>0$ there  exists a solenoidal, $T$-periodic function $\tilde{\bfu} \in W^{\textcolor{black}{1},2}(W^{m,q}),$ $m\in \mathbb{N},$ $q\in[1,\infty],$ such that 
$$\ba{ll}\medskip
\tilde{\bfu}(x,t)=\bfxi(t)+\bfomega\times\bfx\,,\ (t,\bfx)\in[0,T]\times\partial\Omega\,,\\ \medskip
\tilde{\bfu}(x,t)=0\,,\ \mbox{for all $t\in [0,T]$, all $|\bfx|\ge\rho$, and some $\rho>R_*$}\,,\\ 
\| \tilde{\bfu} \|_{W^{2,2}(W^{m,q})}\leq C\,\left(\|\bfxi\|_{W^{2,2}(0,T)}+|\bfomega|\right)
\,,\ea
$$
where $C=C(\Omega,m,q)$.  Moreover
\be
\left|\int_{\Omega_R}\bfv\cdot\nabla\tilde{\bfu}(t)\cdot\bfv\right|\le \varepsilon \|\nabla\bfv\|_2^2\,,\ \ \mbox{for all $\bfv\in H(\Omega_R)\cap W_0^{1,2}(\Omega_R)$}\,.
\eeq{LH}
\EL{ext}
The proof of the following lemma is given in \cite[Lemma 3]{GaSiRo}
\Bl
Let $\bfw \in H(\Omega_{R})\cap W_0^{1,2}(\Omega_R)\cap W^{2,2}(\Omega_{R}).$ Then 
the following properties hold. 
\begin{itemize}
%\item[{\rm (i)}] $\bfomega \times \bfx \cdot \nabla \bfw|_{\partial B_R}  =  0\,;$
\item[{\rm (i)}] $(\bfomega \times \bfw - \bfomega \times \bfx \cdot \nabla \bfw) \in H(\Omega_{R})$\,;  
\item[{\rm (ii)}] 
$\Int{\Omega_R}{} (\bfomega \times \bfw - \bfomega \times \bfx \cdot \nabla \bfw)\cdot P\Delta \bfw\smallskip \\
\hspace*{1.2cm}=  \Int{\partial\Omega}{} \left[  -\bfn\cdot\nabla \bfw \cdot (\bfomega\times \bfx\cdot\nabla \bfw) +\half |\nabla \bfw|^2\bfomega\times \bfx\cdot \bfn\right]-\Int{\Omega_R}{}\nabla(\bfomega\times \bfw):\nabla \bfw \,.$
\end{itemize}
\EL{bello}
\par
We conclude by recalling the next lemma  that ensures  suitable existence and uniqueness properties for a  linear Cauchy problem \cite[Theorem {VIII.4.4}]{Gab}
\Bl
Let $\bfcalg$ be a second-order tensor field in $\mathbb{R}^3 \times (0,\infty)$ such that
$$
[\!]\bfcalg(t)[\!]_{\infty,2,\lambda} + \essup{t\geq 0}\|\nabla \cdot \bfcalg(t)\|_2 < \infty\,,
$$
and let $\bfh\in L^{\infty,q}(\real^3\times (0,\infty))$, $q\in (3,\infty)$, with spatial support contained in a ball of radius $\rho$, some $\rho>0$, centered at the origin.
Then, the problem 
\be
\ba{ll}\ms\left.\ba{rl}\ms
{\bfw}_{t}&=\Delta\bfw+\lambda\,\partial_1\bfw-\nabla\phi+\nabla\cdot\bfcalg+\bfh\\
\nabla\cdot\bfw&=0\ea\right\}\ \ \mbox{in $\real^3\times (0,T)$}\\
\bfw(x,0)=\0\,,
\ea
\eeq{VIII.4.1}
has one and only one solution such that for all $T>0$,
\begin{equation}
\bfw\in L^2(0,T;W^{2,2})\,, \ \bfw_t\in L^2(0,T;L^2)\,;\ \ \nabla\phi\in L^{2}(0,T;L^2).
\label{estimunn1}
\end{equation}
Moreover, 
$$
[\!] \bfw(t) [\!]_{\infty,1,\lambda}
% + \essup{t\geq 0}\| \phi(t) \|_r
<\infty\,,
$$
%for arbitrary $r > \frac{3}{2}$, 
and the following inequality holds:
\begin{equation}
[\!] \bfw(t) [\!]_{\infty,1,\lambda} 
%+ \essup{t\geq 0}\| \phi(t) \|_r 
\leq C\,\left([\!]\bfcalg(t)[\!]_{\infty,2,\lambda} +\essup{t\geq 0}\|\bfh(t)\|_q\right)\, \,  
\label{estimunn2}
\end{equation}
with $C=C(q,\rho,B)$, whenever ${\lambda}\in [0,B]$, for some $B>0$. 
\EL{VIII.4.4}

\setcounter{section}{1}\setcounter{equation}{0}
\section{Linear Problem}
The purpose of this section is to show existence and uniqueness of $T$-periodic solutions, in appropriate function classes, to the following set of linear equations
\be\ba{cc}\smallskip\left.\ba{ll}\medskip
\bfu_t-\bfV(t)\cdot\nabla\bfu+\bfomega\times\bfu=\Delta\bfu-\nabla {p}+\bff\\
\Div\bfu=0\ea\right\}\ \ \mbox{in $\Omega\times (0,T)$}\\
\bfu(x,t)=\bfV(t)\,,\ \ (x,t)\in \partial\Omega\times [0,T]\,,
\ea
\eeq{4.7}
where $\bfV(t):=\bfxi(t)+\bfomega\times\bfx$, with $\bff=\bff(x,t)$,  $\bfxi=\bfxi(t)$  suitably prescribed $T$-periodic functions, and $\bfomega\in\real^3$.\smallskip\par We begin to prove an existence theorem by employing the approach given in \cite{GS1}, that combines Galerkin method with the ``invading domains" technique. Thus, let $\cals=\{R_m, \ m\in\nat\}$ be an increasing,  unbounded sequence of positive numbers with 
$R_1 >R_*$, and denote by  $\{\Omega_{R}, \ R\in\cals\}$   the corresponding sequence of bounded domains with $\cup_{R\in \cals}\Omega_{R}=\Omega$.  
In each $\Omega_R$, $R\in\cals$, we shall look for a $T-$periodic solution $\bfu_R:=\bfv_R+\tilde{\bfu}$ with $\tilde{\bfu}$  given by \lemmref{ext} and $\bfv_R$ solution (in the appropriate functional class) to the following problem
\be\ba{cc}\smallskip\left.\ba{ll}\medskip
(\bfv_R)_t-\bfV(t)\cdot\nabla\bfv_R+\bfomega\times\bfv_R=\Delta\bfv_R-\nabla {p}-\tilde{\bfu}\cdot\nabla\bfv_R-{\bfv}_R\cdot\nabla\tilde{\bfu}+\tilde{\bff}\\
\Div\bfv_R=0\ea\right\}\ \ \mbox{in $\Omega_R\times (0,T)$}\\
\bfv_R(x,t)=\0\,,\ \ (x,t)\in \partial\Omega_R\times [0,T]\,,
\ea
\eeq{1.2}
where 
\be
\tilde{\bff}=\bff +  \Delta \tilde{\bfu} -  \tilde{\bfu}_t + \bfV \cdot \nabla \tilde{\bfu} - \bfomega \times \tilde{\bfu}:=\bff+\bff_c\,.
\eeq{sciu}
For each $R \in \cals,$ we introduce the ortho-normal base of  $H(\Omega_R)$,  $\{\bfw _{Ri}\}_{i\in \mathbb{N}}$, 
constituted by the eigenfunctions of the Stokes problem:
\be
P\Delta \bfw_{Rj} = - \lambda_{Rj} \bfw_{Rj}\,,\ \ \ \bfw_{Rj}\in V(\Omega_R)\cap W^{2,2}(\Omega_R)\,,
\label{stokes}
\ee
and look for an ``approximating" solution to \eqref{1.2} of the form  
$$
\bfv_{Rk}(x,t)= \sum _{i=1}^{k}c_{Rki}(t)\bfw_{Ri} (x)\,,
$$ 
where  the coefficients $\bfc_{Rk}=\{c_{Rk1},\cdots,c_{Rkk}\}$  solve the following system of 
equations   
\be
\dot{c}_{Rkj}= \sum_{i=1}^k A_{ij}(t)c_{Rki} + C_j(t)\,,\qquad j=1,...,k  
\eeq{ivpa}
with
$$\begin{array}{ll} \medskip
A_{ij} := - \ds  (\nabla \bfw_{Ri} , \nabla \bfw_{Rj})_{\Omega_R} - ( \bfomega \times  \bfw_{Ri} , \bfw_{Rj})_{\Omega_R} + (\bfV \cdot \nabla \bfw_{Ri}, \bfw_{Rj})_{\Omega_R}  
\\ \medskip
\hspace*{2.5cm}- (\bfw_{Ri} \cdot \nabla \tilde{\bfu}, \bfw_{Rj})_{\Omega_R} - (\tilde{\bfu} \cdot \nabla \bfw_{Ri}, \bfw_{Rj})_{\Omega_R}\\  
\hspace*{2mm}C_j := \ds (\tilde{\bff} , \bfw_{Rj})_{\Omega_R}
\end{array}
$$
In \cite[Lemmas 3.1, 3.2 and 4.1]{GS1} the following result is proved.\footnote{Actually, in the more general nonlinear context.} 
\Bl Let $\bff=\Div\bfcalf\in L^2(L^2)$ with $\bfcalf\in L^2(L^2)$, $\bfxi\in W^{1,2}(0,T)$ be $T$-periodic. Then, for each $k\in\nat$ problem \eqref{ivpa} has at least one $T$-periodic solution $\bfc_{Rk}=\bfc_{Rk}(t)$. Moreover, the approximating solution $\bfv_{Rk}$ satisfies the following uniform estimates
\be\ba{ll}\medskip
\Sup{t\in[0,T]}\left(\|\bfv_{Rk}(t)\|_6+\|\nabla\bfv_{Rk}(t)\|_2\right)+\|D^2\bfv_{Rk}(t)\|_{L^2(L^2)}\\
\hspace*{4.2cm}
\le C\, \left(\|\bff\|_{L^2(L^2)}+\|\bfcalf\|_{L^2(L^2)}+\|\bfxi\|_{W^{1,2}(0,T)}+|\bfomega|\right)\ea
\eeq{Art}
with $C>0$ independent of $k\in\nat$.
\EL{Ar}
Our next objective is to prove further uniform estimates. Precisely, we have the following.
\Bl Let $\bff=\Div\bfcalf$. Suppose $\bff,\bfcalf\in W^{1,2}(L^2)$ and $\bfxi\in W^{2,2}(0,T)$ are $T$-periodic. Then the corresponding $T$-periodic solution $\bfv_{Rk}$ to \eqref{ivpa} obeys the uniform bound:
\be\ba{ll}\medskip
\Sup{t\in[0,T]}\left(\|(\bfv_{Rk})_t(t)\|_6+\|\nabla(\bfv_{Rk})_t(t)\|_2\right)+\|D^2(\bfv_{Rk})_t\|_{L^2(L^2)}\\ \hspace*{6cm}
\le C\, \left(\|\bff\|_{W^{1,2}(L^2)}+\|\bfcalf\|_{W^{1,2}(L^2)}+\|\bfxi\|_{W^{2,2}(0,T)}+|\bfomega|\right)\,,\ea
\eeq{Art1}
with $C=C(T,V_0)$, where $V_0$ is any fixed upper bound of $\|\bfxi\|_{W^{2,2}(0,T)}+|\bfomega|$.
\EL{Ar1} 
{\em Proof.} We take the time-derivative of both sides of \eqref{ivpa}, multiply the resulting equations by $\dot{ c}_{Rkj}$,  sum over the index $j$ from 1 to $k$, and integrate by parts over $\Omega_R$. We thus get (with $(\cdot,\cdot)\equiv(\cdot,\cdot)_{\Omega_R}$, $\bfv\equiv \bfv_{Rk}$):
\be\ba{rl}\medskip
\ode{}t\|\bfv_t\|_2^2&\!\!\!=-\|\nabla\bfv_t(t)\|_2^2+(\dot{\bfxi}\cdot\nabla\bfv-\tilde{\bfu}_t\cdot\nabla\bfv+\bfv_t\cdot\nabla\tilde{\bfu}+\bfv\cdot\nabla\tilde{\bfu}_t-\bff_c,\partial_t\bfv)-(\bfcalf_t,\nabla\bfv_t)\\
&\!\!\!:=-\|\nabla\bfv_t(t)\|_2^2+\Sum{i=1}6\,\mathscr I_i 
\,.
\ea
\eeq{oo}
By classical embedding theorems we deduce
$$
|\mathscr I_2|+|\mathscr I_4|\le c\,\| \tilde{\bfu}\|_{W^{2,2}(W^{3,2})}\left(\|\nabla\bfv\|_2+\|\bfv\|_{2,K}\right)\|\bfv_t\|_{2,K}\,,
$$
where $K$ is the bounded (spatial) support of $\tilde{\bfu}$, and so, with the help of \lemmref{ext} and Poincar\'e inequality we infer
\be
|\mathscr I_1|+|\mathscr I_2|+|\mathscr I_4|\le c\,\|\nabla\bfv\|_2\|\nabla\bfv_t\|_2\,,
\eeq{o0}
with $c=c(V_0)$.
Likewise,  using \eqref{sciu},  \lemmref{ext} and again Poincar\'e inequality we get
\be
|\mathscr I_5|\le c \,(\|\bfxi\|_{W^{2,2}(0,T)}+|\bfomega|)\|\nabla\bfv_t\|_2\,.
\eeq{o2}
Also, by \eqref{LH} we may choose $\tilde{\bfu}$ such that
\be
|\mathscr I_3|\le \half\|\nabla\bfv_t\|_2^2\,.
\eeq{o2}
Finally, by Schwarz inequality,
\be
|\mathscr I_6|\le \|\bfcalf_t\|_2\,\|\nabla\bfv_t\|_2\,.
\eeq{o3}
Replacing in \eqref{oo} the estimates \eqref{o0}--\eqref{o3}, we find
$$
\ode{}t\|\bfv_t\|_2^2+\half\|\nabla\bfv_t\|_2^2\le 
c\,\big(\|\bfxi\|_{W^{2,2}(0,T)}+|\bfomega|+\|\nabla\bfv\|_2+\|\bfcalf_t\|_2\big)\|\nabla\bfv_t\|_2\,,
$$
which, thanks to \lemmref{Ar} produces
$$
\ode{}t\|\bfv_t\|_2^2+\half\|\nabla\bfv_t\|_2^2\le 
c\,\big(\|\bfxi\|_{W^{2,2}(0,T)}+|\bfomega|+\|\bff\|_{L^2(L^2)}+\|\bfcalf\|_{W^{1,2}(L^2)}\big)\|\nabla\bfv_t\|_2\,,
$$
Integrating both sides of the latter from 0 to $T$ and using the $T$-periodicity of $\bfv$ we easily conclude
\be
\|\nabla\bfv_t\|_{L^2(L^2)} \le 
c\,\big(\|\bfxi\|_{W^{2,2}(0,T)}+|\bfomega|+\|\bff\|_{L^2(L^2)}+\|\bfcalf\|_{W^{1,2}(L^2)}\big)\,,
\eeq{1.11}
where $c>0$ depends on $V_0$ but is independent of $\bfv$. We next take the time derivative of both sides of \eqref{ivpa}, multiply the resulting equations by $-\lambda_{Rj}\dot{ c}_{Rkj}$ sum over the index $j$ from 1 to $k$, and integrate by parts over $\Omega_R$. We thus get 
\be\ba{ll}\medskip
\ode{}t\|\nabla\bfv_t\|_2^2
=-\|P\Delta\bfv_t\|_2^2 +\big(\bfomega\times\bfx\cdot\nabla\bfv_t -\bfomega\times\bfv_t, P\Delta\bfv_t\big)+\big(\dot{\bfxi}\cdot\nabla\bfv+\bfxi\cdot\nabla\bfv_t,P\Delta\bfv_t)\\
\hspace*{2.2cm}-\big(\tilde{\bfu}_t\cdot\nabla\bfv+\tilde{\bfu}\cdot\nabla\bfv_t+\bfv_t\cdot\nabla\tilde{\bfu}-\bfv\cdot\nabla\tilde{\bfu}_t+\tilde{\bff}_t,P\Delta\bfv_t\big)\,.
\ea
\eeq{1.12}
Clearly, by Cauchy-Schwarz, 
\be
|(\dot{\bfxi}\cdot\nabla\bfv+\bfxi\cdot\nabla\bfv_t,P\Delta\bfv_t)|\le c\,\|\bfxi\|_{W^{2,2}(0,T)}^2(\|\nabla\bfv\|_2+\|\nabla\bfv_t\|_2)^2+\mbox{$\frac14$}\|\Delta\bfv_t\|_2^2\,.
\eeq{1.13}
Furthermore, by arguing as done previously, and again by Cauchy-Schwarz, we can show
\be\ba{ll}\medskip
|\big(\tilde{\bfu}_t\cdot\nabla\bfv+\tilde{\bfu}\cdot\nabla\bfv_t+\bfv_t\cdot\nabla\tilde{\bfu}-\bfv\cdot\nabla\tilde{\bfu}_t+\tilde{\bff}_t,P\Delta\bfv_t\big)|\\
\hspace*{3cm}\le 
c\,[(\|\bfxi\|_{W^{2,2}(0,T)}+|\bfomega|+\|\nabla\bfv\|_2+\|\nabla\bfv_t\|_2+\|\bff_t\|_2)]^2+\mbox{$\frac14$}\|\Delta\bfv_t\|_2^2\,,\ea
\eeq{1.14}
with $c=c(V_0)$. 
Also, by \lemmref{bello}(ii) with $\bfw\equiv\bfv_t$, we get
\be
\big|\big(\bfomega\times\bfx\cdot\nabla\bfv_t -\bfomega\times\bfv_t, P\Delta\bfv_t\big)\big|\le c\,|\bfomega|\big( \|\nabla\bfv_t\|^2_2+\int_{\partial\Omega}|\nabla\bfv_t|^2\big)\,.
\eeq{1.15}
We now recall the well-known trace inequality \cite[Theorem II.4.1]{Gab} 
\be
\int_{\partial\Omega}|\nabla\bfv_t|^2\le c_\varepsilon\,\|\nabla\bfv_t\|_2^2+\varepsilon\,\|D^2 \bfv_t|_2^2\,,\ \ \varepsilon>0\,,
\eeq{1.16}
and Heywood inequality \cite[Lemma 1]{Hey}
\be
\|D^2\bfv_t\|_2\le c_0\,\big(\|P\Delta\bfv_t\|_2+\|\nabla\bfv_t\|_2\big)\,,
\eeq{1.17}
where the constant $c_0$ is {\em independent} of $R$.
Thus, if we use \eqref{1.16} and \eqref{1.17} into \eqref{1.15} we get with a suitable choice of $\varepsilon$
\be
\big|\big(\bfomega\times\bfx\cdot\nabla\bfv_t -\bfomega\times\bfv_t, P\Delta\bfv_t\big)\big|\le c\, \|\nabla\bfv_t\|^2_2+\mbox{$\frac14$}\|P\Delta\bfv_t|_2^2\,,
\eeq{1.18}
with $c=c(V_0)$.
If we replace \eqref{1.13}, \eqref{1.14} and \eqref{1.18} into \eqref{1.12}  and use \lemmref{Ar}  we deduce
$$%\ba{ll}\medskip
\ode{}t\|\nabla\bfv_t\|_2^2
+\mbox{$\frac14$}\|P\Delta\bfv_t\|_2^2\\
\le  
c\,\big(\|\bfxi\|_{W^{2,2}(0,T)}+|\bfomega|+\|\nabla\bfv\|_2+\|\nabla\bfv_t\|_2+\|\bff_t\|_2\big)^2\,,
%\ea
$$
which, in turn, with the help of
 \eqref{Art} and \eqref{1.11} furnishes
\be%\ba{ll}\medskip
\ode{}t\|\nabla\bfv_t\|_2^2
+\mbox{$\frac14$}\|P\Delta\bfv_t\|_2^2\\
\le  
c\,\big(\|\bfxi\|_{W^{2,2}(0,T)}+|\bfomega|+\|\bff\|_{W^{1,2}(L^2)}+\|\bfcalf\|_{W^{1,2}(L^2)}\big)^2\,,
%\ea
\eeq{1.20}
with $c=c(V_0)$. We now observe that, by \eqref{1.11}, there is at least one $\bar{t}\in (0,T)$ such that
$$
\|\nabla\bfv_t(\bar{t})\|_{L^2} \le 
c\,\big(\|\bfxi\|_{W^{2,2}(0,T)}+|\bfomega|+\|\bff\|_{L^2(L^2)}+\|\bfcalf\|_{W^{1,2}(L^2)}\big)\,,
$$
and so, integrating \eqref{1.20} between $\bar{t}$ and arbitrary $t>\bar{t}$, and exploiting the $T$-periodicity of $\bfv$, we readily get 
$$
\sup_{t\in[0,T]}\|\nabla\bfv_t({t})\|_{2} +\int_0^T\|P\Delta\bfv\|_2^2\le 
c\,\big(\|\bfxi\|_{W^{2,2}(0,T)}+|\bfomega|+\|\bff\|_{W^{1,2}(L^2)}+\|\bfcalf\|_{W^{1,2}(L^2)}\big)\,.
$$
The lemma then follows from the latter, \eqref{1.17} and the Sobolev inequality $\|w\|_6\le \gamma\,\|\nabla w\|_2$, $w\in W^{1,2}_0(\Omega)$, with $\gamma$ numerical constant.\hfill$\square$\par
With the two previous lemmas in hand, we are now in a position to prove the following result.
\Bl Suppose $\bff$ and $\bfxi$ satisfy the assumptions of \lemmref{Ar1}. Then, there exists at least one $T$-periodic solution $(\bfu,p)$ to \eqref{4.7}, such that \ \mbox{for arbitrarily given $R>2R_*$}
\be
\bfu\in W^{1,\infty}(L^6\cap D^{1,2})\cap W^{1,2}(D^{2,2})\cap L^\infty(L^\infty(\Omega_{\frac R2}))\,,\ \ p\in L^\infty(L^2(\Omega_R))\cap L^2(D^{1,2})\,.
\eeq{1.21}
Moreover, the following estimate holds
\be
\ba{ll}\medskip
\|\bfu\|_{L^\infty(L^\infty(\Omega_{\frac R2}))}+\|\bfu\|_{W^{1,\infty}(L^6)}+\|\nabla\bfu\|_{W^{1,\infty}(L^2)}+\|D^2\bfu\|_{W^{1,2}(L^2)}+\|p\|_{L^\infty(L^2(\Omega_R))}+\|\nabla p\|_{L^{2}(L^2)}\\ \hspace*{4cm}
\le C \left(\|\bff\|_{W^{1,2}(L^2)}+\|\bfcalf\|_{W^{1,2}(L^2)}+\|\bfxi\|_{W^{2,2}(0,T)}+|\bfomega|\right)\,,\ea
\eeq{1.22}
with $C=C(R,\Omega, T,V_0)$, where $V_0$ is any fixed upper bound of $\|\bfxi\|_{W^{2,2}(0,T)}+|\bfomega|$. Moreover, if $\bfxi(t)\not\equiv 0$, then $(\bfu,p)$ is also unique in the class \eqref{1.21}, provided $\bfxi(t)/|\bfxi(t)|$ does not depend on $t\in[0,T]$\,.
\EL{1.3}
{\em Proof.} Once the estimates \eqref{Art} and \eqref{Art1} have been proved for the ``approximating solution", then by following step by step the argument given in \cite[Sections 3 and 4]{GS1} we can show the existence of a solution $(\bfv_R,p_R)$ to \eqref{1.2}, with $\bfu_R:=\bfv_R+\tilde{\bfu}$ in the class \eqref{1.21},  $\nabla p_R\in L^2(D^{1,2})$, and both satisfying \eqref{1.22}. We next let $R\to\infty$ (along a sequence) and, again by entirely following the procedure of \cite{GS1}, we show that $\bfu_R$ tends (in suitable topology) to a solution $\bfu$ to \eqref{4.7} with the regularity properties given in \eqref{1.22}. To such a $\bfu$ one can then associate a pressure field $p\in L^2(D^{1,2})$ with $\bfu$ and $\nabla p$ satisfying \eqref{1.22}. Thus, to complete the proof of the existence property, it remains to show the local bound on $p$ and $\bfu$. To this end, for an arbitrarily fixed $R>2R_*$, let us multiply both sides of \eqref{4.7} by $\bfpsi\in W_0^{1,2}(\Omega_R)$, and integrate by parts as necessary. We thus get
\be
(\bfu_t,\bfpsi)-(\bfV(t)\cdot\nabla\bfu,\bfpsi)+(\bfomega\times\bfu,\bfpsi)+(\nabla\bfu,\nabla\bfpsi)=(p,\Div\bfpsi)+(\bff,\bfpsi)\,.
\eeq{1.24}
We modify $p$ by a ($T$-periodic)  function of time in such a way that $\int_{\Omega_R}p=0$, and, for fixed $t\in[0,T]$ we  choose $\bfpsi$ as solution to \eqref{Bog} with $f\equiv p$. From \eqref{1.24} we thus infer
$$
\|p\|_{2,\Omega_R}^2\le c\,\left(\|\bfu_t\|_{2,\Omega_R}+\|\nabla\bfu\|_{1,2,\Omega_R}+\|\bff\|_{2,\Omega_R}\right)\|\bfpsi\|_{1,2,\Omega_R}
$$
with $c=c(R,V_0)$. Thus, using the properties of $\bfpsi$ given in  \lemmref{1.1_1} along with the estimate \eqref{1.22} for $\bfu$, we prove the desired bound for $p$. Next,  for fixed $t\in[0,T]$, we write \eqref{4.7} as a Stokes problem:
\be
\ba{cc}\smallskip\left.\ba{ll}\medskip
\Delta\bfu-\nabla {p}=\bfF\\
\Div\bfu=0\ea\right\}\ \ \mbox{in $\Omega\times \{t\}$}\\
\bfu(x,t)=\bfV(t)\,,\ \ (x,t)\in \partial\Omega\times \{t\}\,,
\ea
\eeq{stok}
where
$$
\bfF:=\bfu_t-\bfV(t)\cdot\nabla\bfu+\bfomega\times\bfu-\bff\,.
$$
By what we already showed, we have
$$\ba{ll}\medskip
\|\bfF\|_{L^\infty(L^2(\Omega_R))}+\|\bfu\|_{L^\infty(W^{1,2}(\Omega_R))}+\|p\|_{L^\infty(L^2(\Omega_R))}\\
\hspace*{3cm}\le  C \left(\|\bff\|_{W^{1,2}(L^2)}+\|\bfcalf\|_{W^{1,2}(L^2)}+\|\bfxi\|_{W^{2,2}(0,T)}+|\bfomega|\right)
\ea
$$
and so, using the local estimates for problem \eqref{stok} \cite[Theorem IV.5.1]{Gab}:
$$
\|D^2\bfu(t)\|_{2,\Omega_{\frac R2}}\le c\,(\|\bfF(t)\|_2+\|\bfu\|_{1,2,\Omega_R}+\|p\|_{2,\Omega_R}+|\bfxi(t)|+|\bfomega|)\,,
$$
along with the embedding inequality
$$
\|\bfu\|_{\infty,\Omega_{\frac R2}}\le c\,\|\bfu\|_{2,2,\Omega_R}
$$
allows us to obtain the desired estimate for $\|\bfu\|_{L^\infty(L^\infty(\Omega_{\frac R2}})$, thus
concluding the proof of existence. We shall now prove  uniqueness, namely,  that $\bfu\equiv\nabla p\equiv\0$ is the only $T$-periodic solution   in the class  \eqref{1.21} to the following system
\be\ba{cc}\smallskip\left.\ba{ll}\medskip
\bfu_t-\bfV(t)\cdot\nabla\bfu+\bfomega\times\bfu=\Delta\bfu-\nabla {p}\\
\Div\bfu=0\ea\right\}\ \ \mbox{in $\Omega\times (0,T)$}\\
\bfu(x,t)=\0\,,\ \ (x,t)\in \partial\Omega\times [0,T]\,.
\ea
\eeq{47}
To this end, we observe that by formally applying the divergence operator on both sides of \eqref{47}$_1$ and taking into account \eqref{47}$_2$ and \lemmref{bello}(i), we deduce
that $p$ obeys the following Neumann problem for a.a. $t\in [0,T]$ in the distributional sense
\be
\Delta p=0\ \ \mbox{in $\Omega$}\,,\ \ \pde p\bfn=-\curl(\psi\,\curl\bfu)\cdot\bfn\ \ \mbox{at $\partial\Omega$},
\eeq{2.17}
where $\psi$ is a smooth function of bounded support that is 1 in a neighborhood of $\partial\Omega$,
and we used the identity $\Delta\bfu=-\curl\curl\bfu.$ Employing well-known results on the Neumann problem \cite[Theorem III.3.2]{Gab} and the fact that $\bfu$ is in the class \eqref{1.21}, we get
$$
\|\nabla p\|_{L^2(L^q)}\le c\,\|\curl\curl\bfu\|_{L^2(L^q(K))}+\|\curl\bfu\|_{L^2(L^q(K))}\,,\ \ \mbox{all $q\in (1,2]$}\,,
$$
with $K=\supp(\psi)$. 
From this and Sobolev inequality, we then infer (by possibly adding to $p$ a suitable $T$-periodic function of time) that the following property holds
\be
p\in L^2(L^r)\,,\ \ \mbox{all $r\in(3/2,6]$}\,.
\eeq{2.20}
Let $\psi_R=\psi_R(x)$ be the function defined in \lemmref{1.1}.
We dot-multiply both sides of \eqref{47}$_1$ by $\psi_R\bfu$, and integrate by parts over $\Omega\times(0,T)$. Observing that, since $\psi_R=\psi_R(|x|)$, 
$$
\bfomega\times\bfx\cdot\nabla\psi_R=0
$$
and 
that $\bfu\in L^\infty(L^2(\Omega\rho))$, all $\rho\ge R_*$, by using also $T$-periodicity we show 
\be\ba{rl}\medskip
\Int0T\Int\Omega{}\psi_R\,|\nabla\bfu|^2&\!\!\!=-\half\Int0T\Int{\Omega^{\frac{R}{\sqrt2}}}{}\nabla\psi_R\cdot\bfxi(t)|\bfu|^2+\Int0T\Int{\Omega^{\frac{R}{\sqrt2}}}{} p\,\nabla\psi_R\cdot\bfu\\ &\!\!\!\!:= -\half I_{1R}+I_{2R}\,.
\ea\eeq{2.22}
From Schwarz inequality, the properties of $\psi_R$, and  \eqref{1.21} we get
$$
|I_{2R}|\le C_1\sup_{t\in [0,T]}\|\nabla\bfu(t)\|_2\int_0^T\|p(t)\|_{2,\Omega^{\frac{R}{\sqrt2}}}\,, 
$$
which, by \eqref{2.20}, furnishes
\be
\lim_{R\to\infty}|I_{2R}|=0\,.
\eeq{2.23}
Next, if $\bfxi(t)\not\equiv\0$, under the given assumption we take, without loss of generality, $\bfxi(t)=\xi(t)\,\bfe_1$. As a consequence, we get
$$
I_{1R} =\Int{0}T\Int{\Omega_{\frac{R}{\sqrt2}}}{}\xi(t)\,\pde{\psi_R}{x_1}\,|\bfu|^2\,.
$$
By H\"older inequality 
 and the summability properties of $\partial\psi_R/\partial x_1$ we then show
$$
|I_{1R}|\le c\,\|\bfxi\|_{W^{1,2}(0,T)}\, \int_0^T\|{\bfu}\|_{6,\Omega^{\frac{R}{\sqrt2}}}^2\,,
$$
which, 
in view of \eqref{1.21}, implies
\be
\lim_{R\to\infty}|I_{1R}|=0\,.
\eeq{2.24}
Uniqueness then follows by letting $R\to\infty$ in \eqref{2.22} and using \eqref{2.23}--\eqref{2.24}. The lemma is completely proved.
\hfill$\square$\par

\Bl Suppose that, at all times $t\in[0,T]$, the vectors $\bfxi(t)$ and $\bfomega$ are parallel. Let $\bfe_1$ be their common direction,  and orient $\bfe_1$ in such a way that setting   $\lambda\,\bfe_1:=(1/T)\int_0^T\bfxi(t){\rm d}t$, it is $\lambda\ge 0$. Denote by $(\bfu,p)$ the solution to \eqref{4.7} given in \lemmref{1.3}. Then, if, in addition, $[\!]\bfcalf[\!]_{\infty,2,\lambda}<\infty$ it follows that  $[\!]\bfu[\!]_{\infty,1,\lambda}<\infty$, and, moreover,
$$
[\!]\bfu[\!]_{\infty,1,\lambda}\le C\,\big([\!]\bfcalf[\!]_{\infty,2,\lambda}+\|\bff\|_{W^{1,2}(L^2)}+\|\bfcalf\|_{W^{1,2}(L^2)}+\|\bfxi\|_{W^{2,2}(0,T)}+|\bfomega|\big)\,,
$$
where $C=C(\Omega,T,V_0)$, whenever $(\|\bfxi\|_{W^{2,2}(0,T)}+|\bfomega|)\in [0,V_0]$, for some $V_0>0$\,.
\EL{1}
{\em Proof.} Let $\psi$ be the ``cut-off" function introduced in \eqref{2.17}, and  let $\bfz$ be a solution to problem \eqref{Bog} with $f\equiv-\nabla\psi\cdot\bfu$. Since $\int_{\Omega_{\bar{R}}}f=0$, where $\Omega_{\bar{R}}\supset \supp(f)$,  \lemmref{1.1_1} guarantees the existence of such a $\bfz$. Thus, setting
\be
\bfw:=\psi\,\bfu+\bfz\,,\ \ {\sf p}:=\psi\,p\,, \ \ \bfcalh=\psi\bfcalf
\eeq{CTM}
from \eqref{4.7} we deduce that $(\bfw,{\sf p})$ is a $T$-periodic solution to the following problem
\be\left.\ba{ll}\medskip
{\partial}_t\bfw-(\bfxi(t)+\bfomega\times\bfx)\cdot\nabla\bfw+\bfomega\times\bfw=\Delta\bfw-\nabla {\sf p}+\Div\bfcalh+\bfg\\
\Div\bfw=0\ea\right\}\ \ \mbox{in ${\real^3}\times (0,T)$}
\,,
\eeq{1}
where
$$
\bfg:=-\bfz_t+\bfxi(t)\cdot\nabla\bfz+\Delta\bfz-2\nabla\psi\cdot\nabla\bfu+p\,\nabla\psi-\bfxi(t)\cdot\nabla\psi\,\bfu\,.
$$
If we extend $\bfz$ to 0 outside its support, we infer that $\bfg$ is of bounded support. Also with the help of \lemmref{1.1_1} and \lemmref{1.3} (with $R>\max\{2R_*,\bar{R}\}$) we easily deduce
\be\ba{ll}\medskip
\Sup{t\ge 0}\|\bfg(t)\|_2\le c\,(\|\bff\|_{W^{1,2}(L^2)}+\|\bfcalf\|_{W^{1,2}(L^2)}+\|\bfxi\|_{W^{2,2}(0,T)}+|\bfomega|)\,,\\
\Div\bfcalh(t)\in L^\infty(L^2)\,,
\ea
\eeq{g}
where, here and in the rest of the proof, $c$ denotes a constant depending, at most, on $\Omega$, $T$ and $V_0$. Set $\bfomega=\omega\,\bfe_1$ and introduce the new variable $\bfy$ defined by 
\be
\bfy=\bfQ(-t)\cdot(\bfx-\bfx_0(t))
\eeq{2}
where
\be
\bfx_0(t):=\int_0^t(\bfxi(s)-\lambda{\bfe_1})\,{\rm d}s\,,
\eeq{3}
and 
$$
{\bfQ}(t)= \plq \begin{array}{ccc}1&0&0\\0&\cos(\omega t)&\sin(\omega t)
\\0&-\sin(\omega t)&\cos(\omega t) \end{array}\prq\,.
$$
We observe the simple but important property that, since $(1/T)\int_0^T\big(\bfxi(t)-\lambda\,{\bfe_1}\big){\rm d}t=\0$, one can show the existence of a constant $M=M(T,V_0)$ such that (see \cite[Proposition 1]{GARMA}) 
\be
\sup_{t\ge 0}|\bfx_0(t)|\le M\,
\eeq{GD0}
which, by \eqref{2} implies, in particular,
\be
|\bfx|-M\le |\bfy|\le |\bfx|+M\,.
\eeq{GD}
Furthermore, by \eqref{GD0},  \eqref{GD} and the fact that \be\bfQ(t)\cdot\bfe_1=\bfe_1\,, \ \mbox{for all $t\in\real$}\,,
\eeq{Qq} 
we show
\be\ba{rl}\medskip
(1+|x|)(1+2\lambda\,s(x))&\!\!\!\le (1+|y|+M)\big(1+2\lambda\,s(y)+2\lambda\,(M+x_{01}(t))\big)
\\
&\!\!\!\le c\,(1+|y|)\,\big(1+2\lambda\,s(y)\big)\,,
\ea
\eeq{LJ}
and, likewise,
\be
(1+|y|)(1+2\lambda\,s(y))\le c\, (1+|x|)\big(1+2\lambda\,s(x)\big)
\,.
\eeq{VFM}
Setting
\be\ba{ll}\medskip
\bfv(\bfy,t)=\bfQ(-t)\cdot\bfw(\bfQ(t)\cdot\bfy+\bfx_0(t),t),\\ \medskip
{\sf p}(\bfy,t)=p(\bfQ(t)\cdot\bfy+\bfx_0(t),t), 
\\  \medskip
\bfh=\bfQ(-t)\cdot\bfg(\bfQ(t)\cdot\bfy+\bfx_0(t),t)\\
\bfcalg(\bfy,t)=\bfQ(-t)\cdot\bfcalh(\bfQ(t)\cdot\bfy+\bfx_0(t),t)\cdot\bfQ(t)\,,
\ea
\eeq{4}
from \eqref{1} with the help of \eqref{Qq} we easily deduce that $(\bfv,{\sf p})$ is a solution to the following Cauchy problem
\be\ba{cc}\medskip\left.\ba{ll}\medskip
\bfv_t-\lambda\,\partial_1\bfv=\Delta\bfv-\nabla {\sf p}+\Div\bfcalg+\bfh\\
\Div\bfv=0\ea\right\}\ \ \mbox{in ${\real^3}\times (0,\infty)$}
\,,\\
\bfv(x,0)=\bfw(x,0)\,,
\ea
\eeq{5}
where, of course, all spatial differential operators act now on the $y$-variable.
We look for a solution to \eqref{5} of the form $(\bfv_1+\bfv_2, {\sf p}_1+{\sf p}_2)$ where
\be\ba{cc}\medskip\left.\ba{ll}\medskip
(\bfv_1)_t-\lambda\,\partial_1\bfv_1=\Delta\bfv_1-\nabla {\sf p}_1+\Div\bfcalg+\bfh\\
\Div\bfv_1=0\ea\right\}\ \ \mbox{in ${\real^3}\times (0,\infty)$}
\,,\\
\bfv_1(x,0)=\0\,,
\ea
\eeq{6}
and
\be\ba{cc}\medskip\left.\ba{ll}\medskip
(\bfv_2)_t-\lambda\,\partial_1\bfv_2=\Delta\bfv_2-\nabla {\sf p}_2\\
\Div\bfv_2=0\ea\right\}\ \ \mbox{in ${\real^3}\times (0,\infty)$}
\,,\\
\bfv(x,0)=\bfw(x,0)\,.
\ea
\eeq{7}
From $\eqref{4}_{4}$,  \eqref{VFM} and \eqref{CTM}$_3$ we infer
\be 
[\!]\bfH[\!]_{\infty,2,\lambda}\le C\,[\!]\bfF[\!]_{\infty,2,\lambda}
\eeq{1.45}
Moreover, by \eqref{3}$_3$ and \eqref{g}, it follows that 
\be
\Sup{t\ge 0}\|\bfh(t)\|_2\le c\,(\|\bff\|_{W^{1,2}(L^2)}+\|\bfcalf\|_{W^{1,2}(L^2)}+\|\bfxi\|_{W^{2,2}(0,T)}+|\bfomega|)\,.
\eeq{1.46}
As a result, from \lemmref{VIII.4.4} we conclude that \eqref{7}
has one and only one solution such that for all $T>0$,
$$
\bfv_1\in L^2(0,T;W^{2,2})\,, \ (\bfv_1)_t\in L^2(0,T;L^2)\,;\ \ \nabla {\sf p}_1\in L^{2}(0,T;L^2)\,,\ 
[\!] \bfw(t) [\!]_{\infty,1,\lambda}
% + \essup{t\geq 0}\| \phi(t) \|_r
<\infty\,,
$$
satisfying, in addition, the inequality
\be 
[\!] \bfv_1 [\!]_{\infty,1,\lambda} 
%+ \essup{t\geq 0}\| \phi(t) \|_r 
\leq c\,\left([\!]\bfcalf[\!]_{\infty,2,\lambda} +\|\bff\|_{W^{1,2}(L^2)}+\|\bfcalf\|_{W^{1,2}(L^2)}+\|\bfxi\|_{W^{2,2}(0,T)}+|\bfomega|\right)\,. 
\eeq{munn2}
Concerning \eqref{7}, 
since $\bfw(x,0)\in L^6(\real^3)$ (by \eqref{1.21}), it follows that there exists a (unique) solution $(\bfv_2,{\sf p}_2)$ with \cite[Theorem VIII.4.3]{Gab}
\be\ba{ll}\medskip
\bfv_2,\partial_t\bfv_2\, D^2\bfv_2\in L^r([\varepsilon,\tau]\times\real^3)\,,\ \ \mbox{all $\varepsilon\in (0,\tau)$, $\tau>0$, and $r\in[6,\infty]$}\,,\\ 
\|\bfv_2(t)\|_\infty\le C_1\,t^{-\frac14}\|\bfw(0)\|_6\,,\ \ \Sup{t\in(0,\infty)}\|\bfv_2(t)\|_6\le C_1\,\|\bfw(0)\|_6\,.
\ea 
\eeq{15}
In view of the regularity properties of $\bfu$ (and hence of $\bfw$) and those in (\ref{estimunn1}), \eqref{15} for $\bfv_i$, $i=1,2$, respectively, we may use the results proved in \cite[Lemma VIII.4.2]{Gab} to guarantee  $\bfw=\bfv_1+\bfv_2$. As a consequence, due to the $T$-periodicity of $\bfw$ and \eqref{4}$_1$, for  arbitrary positive integer $n$ and $t\in[0,T]$ we obtain 
\be\ba{rl}\medskip
|\bfw(x,t)|(1+|x|)(1+2\lambda\,s(x))&\!\!\!\!=|\bfv(y,t+nT)|(1+|x|)(1+2\lambda\,s(x))\\ &\!\!\!\!\le \big(|\bfv_1(y,t+nT)|+|\bfv_2(y,t+nT)|\big)(1+|x|)(1+2\lambda\,s(x)).\ea
\eeq{17}
Employing  \eqref{LJ}, \eqref{munn2} and \eqref{15}$_2$ in this inequality we get
$$\ba{ll}\medskip
|\bfw(x,t)|(1+|x|)(1+2\lambda\,s(x))
\le c\,\big[(1+|x|)(1+2\lambda\,s(x)) (t+nT)^{-\frac14}\|\bfw(0)\|_6\\
\hspace*{4cm} +[\!]\bfcalf[\!]_{\infty,2,\lambda} +\|\bff\|_{W^{1,2}(L^2)}+\|\bfcalf\|_{W^{1,2}(L^2)}+\|\bfxi\|_{W^{2,2}(0,T)}+|\bfomega|\big]\,\ea
$$ 
so that,  by letting $n\to\infty$ and recalling that, uniformly in $t\ge 0$,  $\bfu(x,t)\equiv \bfw(x,t)$  for $|x|$ sufficiently large ($>\bar{R}$)  we deduce
\be
[\!]\bfu[\!]_{\infty,1,\lambda,\Omega^{\bar{R}}}\le c\,\big([\!]\bfcalf[\!]_{\infty,2,\lambda} +\|\bff\|_{W^{1,2}(L^2)}+\|\bfcalf\|_{W^{1,2}(L^2)}+\|\bfxi\|_{W^{2,2}(0,T)}+|\bfomega|\big)\,. 
\eeq{vfn}
Moreover, by \eqref{1.22}  we have (with $R>2\bar{R}$)
\be
\|\bfu\|_{L^\infty(L^\infty(\Omega_{\frac R2}))}\le c\,\big(\|\bff\|_{W^{1,2}(L^2)}+\|\bfcalf\|_{W^{1,2}(L^2)}+\|\bfxi\|_{W^{2,2}(0,T)}+|\bfomega|\big)
\eeq{vfm}
and 
the desired result then follows from \eqref{vfn} and \eqref{vfm}.
\par\hfill$\square$\par
Combining the  results of \lemmref{1.3} and \lemmref{1} we immediately obtain the following theorem representing the main achievement of this section.
\Bt Let $\bff=\Div\bfcalf$. Suppose $\bff,\bfcalf\in W^{1,2}(L^2)$,   and $\bfxi\in W^{2,2}(0,T)$ are $T$-periodic.
Suppose further that, at all times $t\in[0,T]$, the vectors $\bfxi(t)$ and $\bfomega$ are parallel, and let $\bfe_1$ be their common direction.   We  orient $\bfe_1$ in such a way that setting   $\lambda\,\bfe_1:=(1/T)\int_0^T\bfxi(t){\rm d}t$, it is $\lambda\ge 0$. Then, if also $[\!]\bfcalf[\!]_{\infty,2,\lambda}<\infty$, there exists one and only one $T$-periodic solution $(\bfu,p)$ to \eqref{4.7} such that
\be\ba{ll}\medskip
\bfu\in W^{1,\infty}(L^6\cap D^{1,2})\cap W^{1,2}(D^{2,2})\,,\ \ [\!]\bfu[\!]_{\infty,\lambda,1}<\infty\,,\\ p\in L^\infty(L^2(\Omega_R))\cap L^2(D^{1,2})\,.
\ea
\eeq{1.53}
Moreover, the following estimate holds
$$
\ba{ll}\medskip
\|\bfu\|_{L^\infty(L^\infty(\Omega_{\frac R2}))}+\|\bfu\|_{W^{1,\infty}(L^6)}+\|\nabla\bfu\|_{W^{1,\infty}(L^2)}+\|D^2\bfu\|_{W^{1,2}(L^2)}+[\!]\bfu[\!]_{\infty,1,\lambda}\\ \hspace*{3mm}
+\|p\|_{L^\infty(L^2(\Omega_R))}+\|\nabla p\|_{L^{2}(L^2)}
\le C \left(\|\bff\|_{W^{1,2}(L^2)}\!+\!\|\bfcalf\|_{W^{1,2}(L^2)}\!+\![\!]\bfcalf[\!]_{\infty,2,\lambda}+\|\bfxi\|_{W^{2,2}(0,T)}+|\bfomega|\right),\ea
$$
with $C=C(R,\Omega, T,V_0)$, where $V_0$ is any fixed upper bound of $\|\bfxi\|_{W^{2,2}(0,T)}+|\bfomega|$.
\ET{1.1}
\setcounter{equation}{0}
\section{On the Unique Solvability of the Nonlinear Problem}
The main objective of this section is to study the properties of $T$-periodic solutions to the  full nonlinear problem \eqref{0.1}. 
%\be\ba{cc}\smallskip\left.\ba{ll}\medskip
%{\partial}_t\bfu-\bfxi(t)\cdot\nabla\bfu+\bfu\cdot\nabla\bfu=\Delta\bfu-\nabla {p}+\bfb\\
%\Div\bfu=0\ea\right\}\ \ \mbox{in $\Omega\times (-\infty,\infty)$}\\
%\bfu(x,t)=\bfxi(t)\,,\ \ (x,t)\in \partial\Omega\times (-\infty,\infty)\,,
%\ea
%\eeq{3.1}
%where $\bfb=\bfb(x,t)$ and $\bfxi=\bfxi(t)$ are prescribed $T$-periodic functions. 
This will be achieved by combining the results proved in \theoref{1.1} with a classical contraction mapping argument. To this end, we introduce the Banach space
$$\ba{ll}\medskip
\mathscr S:=\big\{\mbox{$T$-periodic $\bfu$}: \Omega\times [0,T]\mapsto\real^3\,|\\
\qquad\qquad[\!]\bfu[\!]_{\infty,1,\lambda}<\infty\,,\  \bfu\in W^{1,\infty}(L^6\cap D^{1,2})\cap W^{1,2}(D^{2,2})\,;\ \Div\bfu=0\big\}\,, 
\ea$$
endowed with the norm
\be%\ba{ll}\medskip
\|\bfu\|_{\mathscr S}:=[\!]\bfu[\!]_{\infty,1,\lambda}+\|\bfu\|_{W^{1,\infty}(L^6\cap D^{1,2})}+\|\bfu\|_{W^{1,2}(D^{2,2})}
\eeq{3.2}
\Bl Let $\textbf{\textsf{u}},\textbf{\textsf{w}}\in \mathscr S$. Then
$
\textbf{\textsf{u}}\cdot\nabla\textbf{\textsf{w}}\in W^{1,2}(L^2)$ and
$$
\|\textbf{\textsf{u}}\cdot\nabla\textbf{\textsf{w}}\|_{W^{1,2}(L^2)}+\|\textbf{\textsf{u}}\otimes\textbf{\textsf{w}}\|_{W^{1,2}(L^2)}\le c\,\|\textbf{\textsf{u}}\|_{\mathscr S}\|\textbf{\textsf{w}}\|_{\mathscr S}\,.
$$
\EL{3.1}
{\em Proof.} Clearly,
$$
\|\textbf{\textsf{u}}\cdot\nabla\textbf{\textsf{w}}\|_{L^{2}(L^2)}\le [\!]\textbf{\textsf{u}}[\!]_{\infty,1,\lambda}\,\|\nabla\textbf{\textsf{w}}\|_{L^\infty(L^2)}\le \|\textbf{\textsf{u}}\|_{\mathscr S}\|\textbf{\textsf{w}}\|_{\mathscr S}\,,
$$
and, likewise,
$$
\|\textbf{\textsf{u}}\otimes\textbf{\textsf{w}}\|_{L^{2}(L^2)}\le c\,[\!]\textbf{\textsf{u}}[\!]_{\infty,1,\lambda}\,[\!]\textbf{\textsf{w}}[\!]_{\infty,1,\lambda}\le c\,\|\textbf{\textsf{u}}\|_{\mathscr S}\|\textbf{\textsf{w}}\|_{\mathscr S}\,.
$$
Moreover, by using the inequality $\|\nabla w\|_{3}\le c\,\|\nabla w\|_2^{\frac12}\|D^2w\|_2^{\frac12}$ (see \cite[Theorem 2.1]{CrMa}) along with H\"older inequality, we get
$$\ba{rl}\medskip
\|\textbf{\textsf{u}}_t\cdot\nabla\textbf{\textsf{w}}\|_{L^{2}(L^2)}+\|\textbf{\textsf{u}}\cdot\nabla\textbf{\textsf{w}}_t\|_{L^{2}(L^2)}&\!\!\!\le\|\textbf{\textsf{u}}_t\|_{L^{\infty}(L^6)} \|\nabla\textbf{\textsf{w}}\|_{L^2(L^3)}+[\!]\textbf{\textsf{u}}[\!]_{\infty,1,\lambda}\, \|\nabla\textbf{\textsf{w}}_t\|_{L^2(L^2)}\\ \medskip
&\!\!\! \le c\,\big(\|\textbf{\textsf{u}}\|_{\mathscr S}(\|\nabla\textbf{\textsf{w}}\|_{L^\infty(L^{2})}^{\frac12}\|D^2\textbf{\textsf{w}}\|_{L^2(L^{2})}^{\frac12})+\|\textbf{\textsf{u}}\|_{\mathscr S}\|\textbf{\textsf{w}}\|_{\mathscr S}\big)\\
&\!\!\! \le c\,\|\textbf{\textsf{u}}\|_{\mathscr S}\|\textbf{\textsf{w}}\|_{\mathscr S}\,.
\ea
$$
Finally, employing Hardy inequality \cite[Theorem II.6.1]{Gab}
$$\ba{rl}\medskip
\|\textbf{\textsf{u}}_t\cdot\nabla\textbf{\textsf{w}}\|_{L^{2}(L^2)}+\|\textbf{\textsf{u}}\cdot\nabla\textbf{\textsf{w}}_t\|_{L^{2}(L^2)}&\!\!\!\le [\!]\textbf{\textsf{w}}[\!]_{\infty,1,\lambda}\|\textbf{\textsf{u}}_t/|x|\|_{L^2(L^2)}+[\!]\textbf{\textsf{u}}[\!]_{\infty,1,\lambda}\|\textbf{\textsf{w}}_t/|x|\|_{L^2(L^2)}\\ \medskip
&\!\!\!\le c\,\left([\!]\textbf{\textsf{w}}[\!]_{\infty,1,\lambda}\|\nabla\textbf{\textsf{u}}_t\|_{L^2(L^2)}+[\!]\textbf{\textsf{u}}[\!]_{\infty,1,\lambda}\|\nabla\textbf{\textsf{w}}_t\|_{L^2(L^2)}\right)\\
&\!\!\!\le c\,\|\textbf{\textsf{u}}\|_{\mathscr S}\|\textbf{\textsf{w}}\|_{\mathscr S}\,.
\ea
$$
The proof of the lemma is completed.
\par\hfill$\square$\par
We are now in a position to prove the main result of this paper.
\Bt
Let $\bfcalb=\bfcalb(x,t)$, $\bfxi=\bfxi(t)$ be given $T$-periodic functions such that  
$$\bfcalb,\,\ \bfb:=\Div \bfcalb \in W^{1,2}(L^{2})\,, \ \ [\!]\bfcalb[\!]_{\infty,2,\lambda}<\infty\,,\ \ \bfxi\in W^{2,2}(0,T)\,.
$$
Suppose $\bfxi(t)=\xi(t)\,\bfe_1$ and, without loss of generality, let $T^{-1}\int_0^T{\xi}:=\lambda\ge 0$. Furthermore, let $\bfomega\in\real$ with $\bfomega=\omega\,\bfe_1$.
Then, there exists  $\varepsilon_0=\varepsilon_0(\Omega,T)>0$ such that if 
$$
{\sf D}:=\|\bfb\|_{W^{1,2}(L^2)}+\|\bfcalb\|_{W^{1,2}(L^2)}+[\!]\bfcalb[\!]_{\infty,2,\lambda}+\|\bfxi\|_{W^{2,2}(0,T)}+|\bfomega|<\varepsilon_0\,,
$$
problem \eqref{0.1} has
one and only one  $T$-periodic solution $(\bfu,p)\in \mathscr S\times
 L^2(D^{1,2})$\, with  $\|\bfu\|_{\mathscr S}\le c\,{\sf D}$, for some $c=c(\Omega,T)$.
\ET{3.1}
{\em Proof.} We employ the contraction mapping theorem. To this end, define the map
$$
M:\textbf{\textsf{u}}\in\mathscr S\mapsto \bfu\in\mathscr S\,,
$$
with $\bfu$ solving the linear problem
\be\ba{cc}\smallskip\left.\ba{ll}\medskip
\bfu_t-(\bfxi(t)+\bfomega\times\bfx)\cdot\nabla\bfu+\bfomega\times\bfu=\Delta\bfu-\nabla {p}+\textbf{\textsf{u}}\cdot\nabla \textbf{\textsf{u}}+\bfb\\
\Div\bfu=0\ea\right\}\ \ \mbox{in $\Omega\times (0,T)$}\\
\bfu(x,t)=\bfxi(t)+\bfomega\times\bfx\,,\ \ (x,t)\in \partial\Omega\times [0,T]\,,
\ea
\eeq{lin}
Set 
\be\textbf{\textsf{f}}:=\textbf{\textsf{u}}\cdot\nabla\textbf{\textsf{u}}=\Div(\textbf{\textsf{u}}\otimes \textbf{\textsf{u}}):=\Div\textbf{\textsf{F}}\,,
\eeq{C}
where we used the condition $\Div\textbf{\textsf{u}}=0$. In virtue of \lemmref{3.1},  by assumption, and by the obvious inequality
$$
[\!]\textbf{\textsf{F}}[\!]_{\infty,2,\lambda}\le c_1 [\!]\textbf{\textsf{u}}[\!]_{\infty,1,\lambda}^2\,,\ \ \textbf{\textsf{u}}\in\mathscr S\,,
$$
we infer that
$\textbf{\textsf{F}}$,  $\bfb$,  $\bfxi$ and $\bfomega$ satisfy the assumptions of \theoref{1.1}. Therefore, by that theorem we conclude that the map $M$ is well defined and, in particular, that
\be
\|\bfu\|_{\mathscr S}\le c_2\left(\|\textbf{\textsf{u}}\|_{\mathscr S}^2+{\sf D}\right)\,,
\eeq{3.3}
with $c_2=c_2(\Omega,T,V_0)$. If we now take \be\|\textbf{\textsf{u}}\|_{\mathscr S}<\delta\,,\ \  \delta:=4c_2{\sf D}\,,\ \ {\sf D}<\frac{1}{16 c_2^2}\,,
\eeq{BD}
from \eqref{3.3} we deduce $\|\bfu\|_{\mathscr S}<\half\delta$.
Let $\textbf{\textsf{u}}_i\in\mathscr S$ $i=1,2$, and set 
$$
\textbf{\textsf{u}}:=\textbf{\textsf{u}}_1-\textbf{\textsf{u}}_2\,,\ \ \bfu:=M(\textbf{\textsf{u}}_1)-M(\textbf{\textsf{u}}_2)\,.
$$
From \eqref{lin} we then show
\be\ba{cc}\smallskip\left.\ba{ll}\medskip
\bfu_t-(\bfxi(t)+\bfomega\times\bfx)\cdot\nabla\bfu+\bfomega\times\bfu=\Delta\bfu-\nabla {p}+\textbf{\textsf{u}}_1\cdot\nabla \textbf{\textsf{u}}+\textbf{\textsf{u}}\cdot\nabla\textbf{\textsf{u}}_2\\
\Div\bfu=0\ea\right\}\ \ \mbox{in $\Omega\times (0,T)$}\\
\bfu(x,t)=\0\,,\ \ (x,t)\in \partial\Omega\times [0,T]\,.
\ea
\eeq{line}
Proceeding as in the proof of \eqref{3.3} we can show
$$
\|\bfu\|_{\mathscr S}\le c_2\,\left(\|\textbf{\textsf{u}}_1\|_{\mathscr X}+\|\textbf{\textsf{u}}_2\|_{\mathscr S}\right)\|\textbf{\textsf{u}}\|_{\mathscr S}\,.
$$
As a result, if $\|\textbf{\textsf{u}}_i\|_{\mathscr S}<\delta$, $i=1,2$, from the previous inequality we infer
$$
\|\bfu\|_{\mathscr S}< 2c_2\delta\|\textbf{\textsf{u}}\|_{\mathscr X}\,,
$$
and since by \eqref{BD} $2c_2\delta<1/2$, we may conclude that $M$ is a contraction, which, along with \eqref{BD}, completes the proof of the theorem.
\par\hfill$\square$\par

\ed